\documentclass[12pt,a4paper]{amsart}

\usepackage[english]{babel}
\usepackage[T1]{fontenc}
\usepackage[utf8]{inputenc}

\usepackage{amsmath, amsthm, amscd, amsfonts}
\usepackage{amssymb, tikz}
\usepackage{mathtools}

\usepackage{hyperref}

\usepackage{setspace}
\numberwithin{equation}{section}

\def\rmark{\mbox{$\rm\bf\rule{0.06em}{1.45ex}\kern-0.05em R$}}
\def\pmark{\mbox{$\rm\bf\rule{0.06em}{1.45ex}\kern-0.05em P$}}
\def\nmark{\mbox{$\rm\bf\rule{0.06em}{1.45ex}\kern-0.05em N$}}
\def\vdash{\mbox{$\rm\| \kern-0.13em -$}}
\newcommand{\lusim}[1]{\smash{\underset{\raisebox{1.2pt}[0cm][0cm]{$\sim$}}
{{#1}}}}

\usepackage{pb-diagram}

\def\rmark{\mbox{$\rm\bf\rule{0.06em}{1.45ex}\kern-0.05em R$}}
\def\pmark{\mbox{$\rm\bf\rule{0.06em}{1.45ex}\kern-0.05em P$}}
\def\nmark{\mbox{$\rm\bf\rule{0.06em}{1.45ex}\kern-0.05em N$}}
\def\vdash{\mbox{$\rm\| \kern-0.13em -$}}

\title[Two remarks on Merimovich's model]{Two remarks on Merimovich's model of the total failure of GCH}

\author[M.\  Golshani]{Mohammad Golshani}

\address{Mohammad Golshani, School of Mathematics, Institute for Research in Fundamental Sciences (IPM), P.O.\ Box:
19395--5746, Tehran, Iran.}

\email{golshani.m@gmail.com}

\thanks{The  author's research has been supported by a grant from IPM (No. 99030417).}

\date{}

\begin{document}

%\subjclass[2010]{03E05, 03E35, 03E50, 03E55}

 \maketitle

\begin{abstract}
Let $M$ denote the Merimovich's model \cite{merimovich} in which for each infinite cardinal $\lambda, 2^\lambda=\lambda^{+3}$. We show that in $M$ the following hold:
\begin{enumerate}
\item Shelah's strong hypothesis fails at all singular cardinals, indeed,
\[
\forall \lambda  (\lambda \text{~is a singular cardinal~} \Rightarrow pp(\lambda)=\lambda^{+3}).
\]

\item For each singular cardinal $\lambda$ there is an inner model $N$ of $M$ such that $M$ and $N$ have the same bounded subsets of
$\lambda,$ $\lambda$ is a singular cardinal in $N$, $(\lambda^{+i})^N=(\lambda^{+i})^M$, for $i=1,2,3,$ and $N  \models 2^\lambda=\lambda^{+}$. Thus it is possible to add many new fresh subsets to a singular cardinal
$\lambda$ without adding any new bounded subsets to $\lambda.$

\end{enumerate}
\end{abstract}

\section{introduction}
In \cite{merimovich}, Merimovich has constructed a model of ZFC in which the GCH fails in a uniform way, namely, $2^\lambda=\lambda^{+3}$ for all infinite cardinals $\lambda$. His model is of the form  $M=W_\kappa$, the rank initial segment of $W$ at $\kappa$, where $W$ is a forcing extension of the universe $V$, obtained by first doing a preparation forcing and then forcing with the extender based Radin forcing with interleaved collapses.  In this short note we show two interesting properties of this model  by showing that the following hold in $M$:
\begin{enumerate}
\item Shelah's strong hypothesis fails at all singular cardinals, indeed,
\[
\forall \lambda  (\lambda \text{~is a singular cardinal~} \Rightarrow pp(\lambda)=\lambda^{+3}).
\]

\item For each singular cardinal $\lambda$, there is an inner model $N$ of $M$ such that $M$ and $N$ have the same bounded subsets of
$\lambda,$ $\lambda$ is a singular cardinal in $N$, $(\lambda^{+i})^N=(\lambda^{+i})^M$, for $i=1,2,3,$ and $N  \models 2^\lambda=\lambda^{+}$. Thus it is possible to add many new fresh subsets to
a singular cardinal $\lambda$ without adding any new bounded subsets to $\lambda.$

\end{enumerate}
We assume familiarity with Merimovich's construction and use  the results from \cite{merimovich} without any mention. We also freely use the results or methods used in other references given at the end of the paper.
\section{Merimovich's model}
In this section, we sketch how the model $M$ from the introduction in constructed,  by stating the strategy used to build $M$ and refer to
\cite{merimovich} for the details. Start with a $(\kappa+4)$-strong cardinal $\kappa.$ First we
do a reverse Easton iteration
\[
\langle \langle \mathbb{P}_\alpha: \alpha \leq \kappa+1 \rangle, \langle \lusim{\mathbb{Q}}_\alpha: \alpha \leq \kappa                  \rangle
\]
where at stage $\alpha$, we force with the trivial forcing, except $\alpha$ is an inaccessible cardinal, in which case,
\[
\mathbb{Q}_\alpha=Add(\alpha^+, \alpha^{+4}) \times Add(\alpha^{++}, \alpha^{+5}) \times Add(\alpha^{+3}, \alpha^{+6}).
\]
 Let $\langle \langle G_\alpha: \alpha \leq \kappa+1 \rangle, \langle H_\alpha: \alpha \leq \kappa              \rangle$
be the corresponding generic filter sequence. Now working in $V[G_{\kappa+1}]$, one forces with the extender based Radin forcing with interleaved
collapses $\mathbb{P}_{\bar{E}}$, for a well-chosen coherent sequence of extenders $\bar{E}$. Let $G_{\bar{E}}$ be the corresponding generic filter over $V[G_{\kappa+1}]$. Then $W=V[G_{\kappa+1}][G_{\bar{E}}]$
is a model of ZFC in which $\kappa$ remains inaccessible and for all $\lambda<\kappa,~2^\lambda=\lambda^{+3}$. Then $M=W_\kappa$ is the resulting model of Merimovich.

\section{Shelah's strong hypothesis  fails everywhere in $M$}
In this section we study Shelah's strong hypothesis and show that
\[
M \models ~\forall \lambda  (\lambda \text{~is a singular cardinal~} \Rightarrow pp(\lambda)=\lambda^{+3}).
\]
Thus suppose that $\lambda$ is a singular cardinal in $M$. If $cf^M(\lambda) > \aleph_0$, then by a result of Shelah \cite{shelah}
\[
pp(\lambda)=2^\lambda=\lambda^{+3},
\]
and we are done. If $cf^M(\lambda)=\aleph_0$, then by \cite{gitik}, there is a canonical $\omega$-sequence $\langle \lambda_n: n<\omega \rangle$, added by the Radin club, which is cofinal in $\lambda$\footnote{In this case, either the coherent sequence of extenders has length of cofinality $\omega$, which naturally adds an $\omega$-sequence or its length has cofinality  $\kappa,$ in which case \cite[Lemma 5.13]{gitik} applies.} It is not difficult to show that the sequence $\langle \lambda^{+3}_n: n<\omega \rangle$ admits a scale of length $\lambda^{+3}$ (see also \cite[Lemma 4.10 and the remarks after it]{gitik} for the case of extender based Prikry forcing), and hence $pp(\lambda)=\lambda^{+3}$ holds again in $M$.

\section{Adding many new fresh sets to singular cardinals}
Suppose $\lambda$ is a singular cardinal in $M$. We show that there exists an inner model $N$ of $M$ such that:
\begin{itemize}
\item $M$ and $N$ have the same bounded subsets of
$\lambda,$
\item $\lambda$ is a singular cardinal in $N$,
\item $(\lambda^{+i})^N=(\lambda^{+i})^M$, for $i=1,2,3,$
\item $N  \models 2^\lambda=\lambda^{+}$.
\end{itemize}
Let $C=\langle \kappa_\xi: \xi < \kappa        \rangle$ be the Radin club added by $G_{\bar{E}}$ and let $\xi < \kappa$ be such that
$\lambda=\kappa_\xi.$ Pick some $p =p_l^{\frown}\cdots  ^{\frown} p_k ^{\frown}\cdots ^{\frown} p_0 \in G_{\bar{E}}$ and $\bar{\epsilon}$ such that $\kappa^0(\bar{\epsilon})=\kappa_\xi, p_{l, ..,k} \in \mathbb{P}_{\bar{\epsilon}}$   and
\[
\mathbb{P}_{\bar{E}}/p \simeq \mathbb{P}_{\bar{\epsilon}} / p_{l,...,k} \times \mathbb{P}_{\bar{E}} / p_{k+1, ...,  0}.
\]
Let us first consider the generic extension $V[G_{\kappa_\xi}+1][G_{\bar{\epsilon}}]$ of $V$\footnote{This is well-defined, as the tail forcing $\mathbb{P}_{\bar{E}} / p_{k+1, ...,  0}$ doesn't add new subsets to $\kappa_\xi^{+3}$, so $\mathbb{P}_{\bar{\epsilon}} $ is computed the same in $V[G_{\kappa+1}]$ and $V[G_{\kappa_\xi+1}]$}. Now, working in $V[G_{\kappa_\xi+1}]$,  by \cite{omer}, we can pick some elementary submodel
$A$ of the large part of the universe, so that $A \supseteq V_{\kappa_\xi}$ has size $\kappa_\xi$, $A$ contains all relevant information and $\mathbb{P}_{\bar{\epsilon}} \cap A \lessdot \mathbb{P}_{\bar{\epsilon}}.$ Let
$N$ be the rank initial segment of $V[G_{\kappa_\xi+1}][G_{\bar{\epsilon}} \cap A]$ at $\kappa$. Then $N$ is as required.

\end{document}